\documentclass[12pt,a4paper,twoside]{article}
\usepackage[english]{babel}
\usepackage{amssymb}
\newtheorem{thm}{Theorem}[section]
\newtheorem{cor}[thm]{Corollary}
\newtheorem{lem}[thm]{Lemma}
\newtheorem{prop}[thm]{Proposition}
\newtheorem{defn}[thm]{Definition}
\newtheorem{rem}[thm]{Remark}

\def\b{{\beta}}

\begin{document}

\author{S. Albeverio $^{1},$ Sh. A. Ayupov $^{2},$ \ \ B. A.
Omirov  $^3,$ A.Kh. Khudoyberdiyev$^4$}

\title{\bf n-Dimensional filiform Leibniz algebras of length
(n-1) and their derivations.}

\maketitle

\begin{abstract}
In this work $n$-dimensional filiform Leibniz algebras admitting a
gradation of length $(n-1)$ are classified. Derivations of such
algebras are also described.

\end{abstract}

\medskip
$^1$ Institut f\"{u}r Angewandte Mathematik, Universit\"{a}t Bonn,
Wegelerstr. 6, D-53115 Bonn (Germany); SFB 611, BiBoS; CERFIM
(Locarno); Acc. Arch. (USI),  e-mail: \emph{albeverio@uni-bonn.de}

$^2$ Institute of Mathematics, Uzbekistan Academy of Science, F.
Hodjaev str. 29, 100125, Tashkent (Uzbekistan), e-mail:
\emph{sh\_ayupov@mail.ru, e\_ayupov@hotmail.com, mathinst@uzsci.net}

 $^{3}$ Institute of Mathematics, Uzbekistan
Academy of Science, F. Hodjaev str. 29, 100125, Tashkent
(Uzbekistan), e-mail: \emph{omirovb@mail.ru}

 $^{4}$ Department of Mathematics, National University of Uzbekistan,
Vuzgorogok, 27, 100174, Tashkent (Uzbekistan), e-mail:
\emph{abror\_alg\_85@mail.ru}



\section{Introduction}

The well-known natural gradations of nilpotent Lie and Leibniz
algebras are very helpful in investigations of properties of those
algebras in the general case without restriction on the gradation.
This technique provides a rather deep information on the algebra and
it is more effective when the length of the natural gradation is
sufficiently large. A similar approach was considered in \cite{AO1},
\cite{GK}, \cite{V} (and some other papers). The idea of
consideration of more convenient gradations firstly was suggested in
\cite{GJR1}, \cite{GJR2}. In particular, in the mentioned papers the
authors considered a gradation which has the maximal possible number
of non zero subspaces, i.e. the length of the gradation coincides
with the dimension of the algebra. Actually, gradations with a large
number of non zero subspaces enable us to describe the
multiplication on the algebra more exactly. Indeed, we can always
choose a homogeneous basis and thus the gradation allows to obtain
more explicit conditions for the structural constants.

Moreover, such gradations are useful for the investigation of
cohomologies for the considered algebras because they induce a
corresponding gradation of the group of cohomologies. A similar
approach was considered in \cite{CGO}, \cite{GK}, \cite{O}.

In the present paper we consider gradations which have the length
$(n-1)$ for the classification of complex filiform Leibniz algebras
and apply the obtained results for the description of the algebras
of their derivations.

Recall that an algebra $L$ over a field $F$ is called a {\it Leibniz
algebra} if it satisfies the following Leibniz identity:
$$[x,[y,z]]=[[x,y],z]-[[x,z],y],$$
where $[ \ , \ ]$ denotes the multiplication in $L$.

It is not difficult to check that the variety of Leibniz algebras is
a "non-antisymmetric" generalization of the Lie algebras' variety.

\section{Preliminaries}

Let $L$ be an arbitrary Leibniz algebra of dimension $n$ and let
$\{e_1,...,e_n\}$  be a basis of the algebra $L.$ Then the
multiplication on the algebra $L$ is defined by the products of the
basic elements, namely, as
$[e_i,e_j]=\sum\limits_{k=1}^{n}\gamma_{ij}^ke_k$ , where
$\gamma_{ij}^k$ are the structural constants. Thus the problem of
classification of algebras can be reduced to the problem of finding
a description of the structural constants up to a non-degenerate
basis transformation.

From the Leibniz identity we have the following polynomial
equalities for the structural constants:
$$\sum_{l=1}^{n}(\gamma_{jk}^l\gamma_{il}^m -\gamma_{ij}^l\gamma_{lk}^m +\gamma_{ik}^l\gamma_{lj}^m )=0.$$

But the straightforward description of structural constants is
somewhat cumbersome and therefore usually one has to apply different
methods of investigation.

Let $L$ be a $Z$-graded Leibniz algebra with a finite number of non
zero subspaces, i.e. $L=\bigoplus \limits_{i\in Z} V_i,$  where
$[V_i,V_j]\subseteq V_{i+j}$ for any $i, j \in Z.$

We say that a nilpotent Leibniz algebra $L$ admits a {\it connected
gradation} if $L=V_{k_{1}}\oplus V_{k_{2}}\oplus \dots \oplus
V_{k_{t}},$ where each $V_{i}$ is non-zero for $k_1 \leq i \leq
k_t.$

The number of subspaces $l(\oplus L)=k_t-k_1+1$ is called the length
of the gradation. In the case where $l(\oplus L)=dim L,$ the
gradation is called the maximum length gradation \cite{GJR1},
\cite{GJR2}. If $l(\oplus L)=dim L-1,$ then we have a gradation of
 length $(n-1).$
\begin{defn} The length $l(L)$ of a Leibniz algebra L is defined
as $$l(L) = \max \{l( \oplus L) = k_t - k_1 + 1 \mid L = V_{k_{1}}
\oplus V_{k_{2}} \oplus \dots \oplus V_{k_{t}} \quad \mbox {is a
connected gradation.}\}$$
\end{defn}

This definition means that $l(L)$ is the greatest number of
subspaces from connected gradations which exist in $L.$ Thus, every
Leibniz algebra $L$ has the length at least 1, because we can put
$L=V_0.$

Given an arbitrary Leibniz algebra $L$ we define the lower central
series:
$$L^1=L, \ L^{k+1}=[L^{k},L], \ k\geq 1.$$
\begin{defn} A Leibniz algebra $L$ is said to be filiform if
$dim L^i(L)=n-i,$ where $n=dim L$ and $2\leq i \leq n.$
\end{defn}

\begin{defn} Given a filiform Leibniz algebra $L,$ put
$L_i=L^i/L^{i+1}, \ 1 \leq i\leq n-1,$ and $gr L = L_1 \oplus
L_2\oplus\dots L_{n-1}.$ Then $[L_i,L_j]\subseteq L_{i+j}$ and we
obtain the graded algebra $gr L$. If $gr L$ and $L$ are isomorphic,
denoted by $gr L=L,$ we say that the algebra $L$ is naturally
graded.
\end{defn}

In the following theorem we summarize the results of the works
 \cite{AO1}, \cite{V}.

\begin{thm} \label{thm2.4} Any complex $n$-dimensional naturally graded filiform
Leibniz algebra is isomorphic to one of the following pairwise non
isomorphic algebras:

$NGF_1=\left\{\begin{array}{ll}
[e_1,e_1]=e_{3},&  \\[1mm]
[e_i,e_1]=e_{i+1}, & \  2\leq i \leq {n-1}
\end{array} \right.$ \\[1mm]

$NGF_2=\left\{\begin{array}{ll}
[e_1,e_1]=e_{3}, &  \\[1mm]
[e_i,e_1]=e_{i+1}, & \  3\leq i \leq {n-1}\\[1mm]
\end{array} \right.$ \\

$NGF_3=\left\{\begin{array}{lll} [e_i,e_1]=-[e_1,e_i]=e_{i+1}, &
2\leq i \leq {n-1}\\[1mm]
[e_i,e_{n+1-i}]=-[e_{n+1-i},e_i]=\alpha (-1)^{i+1}e_n & 2\leq i\leq
n-1.
\end{array} \right.$ \\
where all omitted products are equal to zero, $\alpha\in\{0,1\}$ for
even $n$ and $\alpha=0$ for odd $n.$
\end{thm}
As we can see the natural gradations in these algebras have length
equal to $(n-1)$.

Consider the algebra $NGF_1$ and the gradation $V_{-1}=<e_1-e_2>,
 \ V_{n-2}=<e_2>, \ V_i=<e_{n-i}>, \ 0 \leq i \leq n-3.$ Then
$NGF_1=V_{-1}\oplus V_0\oplus\dots \oplus V_{n-2},$ i.e. this
algebra has the maximal length.

It should be noted (see \cite{CGO}) that the algebra $NGF_1$ is the
only filiform non Lie Leibniz algebra of maximal length.

Summarizing the results from \cite{AO1}, \cite{GO} we obtain the
decomposition of all complex filiform Leibniz algebras into three
disjoint classes.

\begin{thm} \label{th2} Any complex $n-$dimensional filiform
Leibniz algebra admits a basis $\{e_1, e_2, \dots, e_n\}$ such that
the table of multiplication of the algebra have one of the following
forms:

$F_1=\left\{\begin{array}{ll}
[e_1,e_1]=e_{3},&  \\[1mm]
[e_i,e_1]=e_{i+1}, & \  2\leq i \leq {n-1}\\[1mm]
[e_1,e_2]=\alpha_4e_4 + \alpha_5e_5+...+ \alpha_{n-1}e_{n-1}+ \theta e_n, & \\[1mm]
[e_j,e_2]=\alpha_4e_{j+2} + \alpha_5e_{j+3}+...+ \alpha_{n+2-j}e_n,
& \ 2\leq j \leq {n-2}
\end{array} \right.$ \\[1mm]
(omitted products are equal to zero)

$F_2=\left\{\begin{array}{ll}
[e_1,e_1]=e_{3}, &  \\[1mm]
[e_i,e_1]=e_{i+1}, & \  3\leq i \leq {n-1}\\[1mm]
[e_1,e_2]=\beta_3e_4 + \beta_4e_5+...+ \beta_{n-1}e_n, & \\[1mm]
[e_2,e_2]=\gamma e_n, & \\[1mm]
[e_j,e_2]=\beta_3e_{j+2} + \beta_4e_{j+3}+...+ \beta_{n+1-j}e_n, & \
3\leq j \leq {n-2}
\end{array} \right.$ \\
(omitted products are equal to zero)

$F_3=\left\{\begin{array}{lll} [e_i,e_1]=e_{i+1}, &
2\leq i \leq {n-1}\\[1mm]
[e_1,e_i]=-e_{i+1}, & 3\leq i \leq {n-1} \\[1mm]
[e_1,e_1]=\theta_1e_n, &   \\[1mm]
[e_1,e_2]=-e_3+\theta_2e_n, & \\[1mm]
[e_2,e_2]=\theta_3e_n, &  \\[1mm]
[e_i,e_j]=-[e_j,e_i] \in lin<e_{i+j+1}, e_{i+j+2}, \dots , e_n>, &
2\leq i \leq n-2,\\[1mm]
& 2 \leq j \leq {n-i}\\[1mm]
[e_i,e_{n+1-i}]=-[e_{n+1-i},e_i]=\alpha (-1)^{i+1}e_n & 2\leq i\leq
n-1.
\end{array} \right.$ \\
where all omitted products are equal to zero, $\alpha\in\{0,1\}$ for
even $n$ and $\alpha=0$ for odd $n.$
\end{thm}

\begin{defn}
A linear transformation $d$ of a Leibniz algebra $L$ is called a
derivation if for any $x, y\in L$
$$d([x,y])=[d(x),y]+[x, d(y)].$$
\end{defn}

The space of all derivations of the algebra $L$ equipped with the
multiplication defined as the commutator, forms a Lie algebra which
is denoted by $Der(L).$

It is clear that the operator of right multiplication $R_x$ by an
 element $x$ of the algebra $L$ (i.e. $R_x(y)=[y,x]$) is also a
derivation. Derivations of this type are called inner derivations.
Similar to the Lie algebras case the set of the inner derivation
$Inn(L)$ forms an ideal of the algebra $Der(L).$

\section{The main result.}
\vspace{1cm}


\begin{prop} \label{prop1} Let $L$ be any filiform Leibniz algebra
from the family $F_1$ of theorem \ref{th2}. If $L$ admits a
gradation with the length $(n-1),$ then this gradation coincides
with the natural gradation. In particular, $L$ is isomorphic to the
algebra $NGF_1.$
\end{prop}

\noindent {\bf Proof.} Let $L$ admit a gradation of length $(n-1)$ :
$L=V_{k_{1}}\oplus V_{k_{1}+1}\oplus \dots \oplus V_{k_{1}+n-2}.$
Take a homogeneous basis $\{x_0, \ x_1, \ \dots, x_{n-1}\},$ where
$x_i\in V_{k_{1}+i-1} \ 1 \leq i \leq n-1,\ x_0\in V_{k_{1}+j-1},$
for some $\ j \ (1 \leq j \leq n-1). $

Since the algebra $L$ is generated by two elements we may suppose
that for some $s$ and $t$ the elements  $x_s, \ x_t$ from the basis
$\{x_0, x_1, x_2, \dots, x_{n-1}\}$ are the generators of $L.$
Expressing these basic elements via the initial basis $\{e_1, e_2,
\dots e_n\}$ we obtain

$$x_s=\sum\limits_{i=1}^{n}a_ie_i, \
x_t=\sum\limits_{j=1}^{n}b_ie_i,$$ where $a_1b_2-a_2b_1\neq 0.$

Without restriction of generality we may assume that $a_1=b_2=1.$

Using the multiplication of the algebra from theorem \ref{th2} we
consider the products:

$[\dots [[\underbrace{x_s,x_s],x_s],\dots
x_s}\limits_{i-times}]=(1+a_2)e_{i+1}+(*)e_{i+2}+ \dots + (*)e_n, \
2 \leq i \leq n-1.$

\textbf{Case 1.} Let $1+a_2\neq 0.$  Set $y_i=[\dots
[[\underbrace{x_s,x_s],x_s],\dots x_s}\limits_{i-times}], \ 1 \leq i
\leq n-1$, $y_n=x_t.$ Then we have $<y_i>\subseteq V_{ik_{s}} \ 1
\leq k \leq n-1$ and from the connectedness of the gradation we
obtain $k_s=\pm 1.$

In the case where $k_s=1$ we have $L=V_1\oplus V_2\oplus \dots
\oplus V_{n-1}.$ The case $k_s=-1,$ i.e. $L=V_{-n+1}\oplus
V_{-n+2}\oplus \dots \oplus V_{-1},$ can be reduced to the case
$k_s=1$ by putting $V_{j}=V_{-j}, \ -1 \leq j \leq n-1.$

Thus, we may assume that $k_s=1.$ Since $l(\oplus L)=n-1,$ we have
the existence of some $p$ ($1\leq p \leq n-1$) such that $y_n\in
V_p.$

Consider the products:
$$[x_t, x_s]=(1+b_1)e_3+(*)e_4 + \dots + (*)e_n,$$
$$[x_s, x_t]=b_1(1+a_2)e_3+(*)e_4 + \dots + (*)e_n.$$
where the asterisks $(*)$ denote appropriate coefficients at the
basic elements.

It is evident that the coefficient at $e_3$ in the products $[x_t,
x_s]$ and $[x_s, x_t]$ are not equal to zero simultaneously
(otherwise $1-a_2b_1=0$). Therefore we have that either $[x_t,
x_s]\neq0$ or $[x_s, x_t]\neq0$ and these elements belong to $V_2,$
i.e. $p=1.$

Thus, we obtain $L=V_1\oplus V_2\oplus \dots \oplus V_{n-1},$ where
$V_1=<y_1, y_n>, \ V_j=<y_j> \ {for} \ 2\leq j \leq n-1,$ but this
gradation coincides with the natural gradation. Therefore in this
case we obtain the algebra $NGF_1.$

\textbf{Case 2.} Let $1+a_2=0.$ Consider the multiplication

$$[x_t, x_t]=b_1(1+b_1)e_3+(*)e_4+\dots + (*)e_n.$$

Since $1-a_2b_1\neq 0,$ one has $b_1\neq -1.$ If $b_1\neq0,$ then
similarly to the case 1 we also obtain the algebra $NGF_1$.
Therefore $b_1=0,$ i.e. $x_t=e_2+b_3e_2+ \dots + b_ne_n.$

Consider the products:

$$[\dots [x_t, \underbrace{x_s], \dots ,
x_s}\limits_{i-times}]=e_{i+2}+(*)e_{i+3}+\dots + (*)e_n, \ 1 \leq i
\leq n-2.$$

Put
$$y_1=x_s, \ y_2=x_t, \ y_i=[\dots, [x_t,
\underbrace{x_s], \dots , x_s}\limits_{{(i-2)}-times}], \ 3 \leq i
\leq n.$$

Then we have $<y_1>\subseteq V_{k_s},$ $<y_2>\subseteq V_{k_t},$
$<y_i>\subseteq V_{k_t+(i-2)k_s}, \ 3 \leq i \leq n.$

From the property of the gradation we have $k_s=\pm 1.$

Let us consider the case where $k_s=1.$ Since $l(\oplus L)=(n-1),$
there exists $p$ ($0\leq p \leq n-2$) such that $y_1\in V_{k_t+p},$
i.e. $k_t=1-p.$

Thus, we obtain $L=V_{1-p}\oplus V_{2-p}\oplus \dots \oplus
V_{n-1-p},$ where $V_{1-p}=<y_2>, \ V_{2-p}= \\ <y_3>,  \dots, \
V_0=<y_{p+1}>, \ V_{1}=<y_1, y_{p+2}>, \ V_{2}=<y_{p+3}>, \dots , \
V_{n-1-p}=\\ <y_n>.$

Consider the product

$[y_1, y_2]=\alpha y_3=\alpha(e_3 +(*)e_4+ \dots + (*)e_n)=\alpha
e_3 +\alpha(*)e_4+ \dots + \alpha(*)e_n.$

On the other hand
$$[y_1, y_2]=(**)e_4+(**)e_5+ \dots + (**)e_n.$$

Comparing the coefficients at $e_3$ we obtain $\alpha=0,$ i.e.
$[y_1, y_2]=0.$ In a similar way we obtain $[y_2, y_2]=0.$

From the equality
$$[y_{i+1}, y_2]=[y_i, [y_1, y_2]]+[[y_i, y_2], y_1]$$
by induction we obtain $[y_i, y_2]=0, \ 3\leq i \leq n-2.$

Thus, we have an algebra with the following table of multiplication:

$[y_i,y_1]=y_{i+1}, \ 2 \leq i\leq n-1, \ [y_1,y_1]=\beta y_{p+3}, \
[y_i,y_j]=0, \ 1 \leq i\leq n, \ 2 \leq j\leq n.$

Applying if necessary a change of basis
$$y_1^{'}=y_1-\beta y_{p+2}, \ y_i^{'}=y_i, \ 2\leq i \leq n$$
we may assume that $\beta=0.$

Without loss of generality we may suppose that $p=0,$ i.e.
$V_1=<y_1, y_2>, V_j=<y_{j+1}>, \ 2 \leq j \leq n-1.$ So, in the
case where $k_s=1$ the gradation of the algebra is the natural
gradation, i.e. $L=NGF_1.$

In the case where $k_s=-1$ by similar arguments we obtain the
natural gradation of the algebra $L$, i.e. the algebra $NGF_1.$
\hfill{\fbox{}}

\begin{defn}
The sets $\mathfrak{R}(L)=\lbrace x\in L \mid [L,x]=0\rbrace$ and
$\mathfrak{L}(L)=\lbrace x\in L \mid [x,L]=0\rbrace$ are called
respectively the right and left annihilators of the algebra $L.$
\end{defn}

It is not difficult to check that the right annihilator is a two-sides ideal of the Leibniz algebra $L.$

Let $L$ be an $n$-dimensional filiform Leibniz algebra and suppose
that $\lbrace e_1, e_2, \dots, e_n\rbrace$ is a basis of $L.$

Further we will need the following lemma.

\begin{lem} $\cite{GO}.$ \label{lem1}
For any $0\leq p \leq n-k, \ 3 \leq k \leq n$ the equality
$$
\sum\limits_{i=k}^{n}a(i)\sum\limits_{j=i+p}^{n}b(i,j)e_j=\sum\limits_{j=k+p}^{n}
\sum\limits_{i=k}^{j-p}a(i)b(i,j)e_j.
$$
holds.
\end{lem}

\begin{prop} \label{prop10} Let $L$ be any filiform Leibniz
algebra from the family $F_2$ of theorem \ref{th2}. If $L$ admits a
gradation with the length $(n-1),$ then it is isomorphic to one of
the following algebras:

$for \ any \ arbitrary \ value \ of \ n$

$NGF_2=\left\{\begin{array}{ll}
[y_1,y_1]=y_{3}, &  \\[1mm]
[y_i,y_1]=y_{i+1}, & \  3\leq i \leq {n-1}\\[1mm]
\end{array} \right.$ \\

$M_1(k)=\left\{\begin{array}{ll}
[y_i,y_1]=y_{i+1}, &  \  1\leq i \leq n-2, \\[1mm]
[y_i,y_n]=y_{k+i-1}, & \  1\leq i \leq n-k, \  3\leq k \leq n-1.
\end{array} \right.$

$for \ any \ odd \ value \ of \ n$

$M_2=\left\{\begin{array}{ll}
[y_i,y_1]=y_{i+1}, &  \  1\leq i \leq n-2, \\[1mm]
[y_i,y_n]=y_{\frac{n+1}{2}+i-1}, & \  1\leq i \leq \frac{n-1}{2},\\[1mm]
[y_n,y_n]=\alpha y_{n-1}, & \alpha\neq0.

\end{array} \right.$

$M_3=\left\{\begin{array}{ll}
[y_i,y_1]=y_{i+1}, &  \  1\leq i \leq n-2 \\[1mm]
[y_n,y_n]=y_{n-1}, &
\end{array} \right.$

where all omitted products are equal to zero and $\{y_1,y_2, \
\dots, \ y_n\}$ is a basis of the corresponding algebra.
\end{prop}

\noindent {\bf Proof.} Let $L$ be an algebra satisfying the
conditions of the proposition. Similar to proposition \ref{prop1} we
can choose homogeneous generators as

$x_s=e_1+a_2e_2+ \dots + a_ne_n$ and $x_t=b_1e_1+e_2+b_3e_3+ \dots
+b_ne_n.$

Consider the multiplications:

$$[\dots\underbrace{[x_s,x_s], \dots,
x_s]}\limits_{i-times}=e_{i+1}+(*)e_{i+2}+\dots+(*)e_n, \ 2\leq i
\leq n-1.$$

Put
$$y_i=[\dots\underbrace{[x_s,x_s], \dots, x_s]}\limits_{i-times},
\ 1 \leq i \leq n-1, \ y_n=x_t.$$ Then we have $<y_i>\subseteq
V_{ik_s} \ 1 \leq i \leq n-1.$ From the connectedness of the
gradation we have that $k_s=\pm 1.$

By the same arguments as above without loss of generality we may
take $k_s=1.$

Since $l(\oplus L)=n-1$ there exists $p$ ($1\leq p\leq n-1$) such
that $y_n\in V_p.$

Consider the products

$[y_n,y_1]=b_1e_3+\sum\limits_{i=4}^{n}b_{i-1}e_i+a_2b_1
\sum\limits_{i=4}^{n}\beta_{i-1}e_i+a_2\gamma
e_n+a_2\sum\limits_{i=3}^{n-2}b_i\sum\limits_{j=i+2}^{n}
\beta_{j-i+1}e_j,$

$[y_1,y_n]=b_1e_3+b_1\sum\limits_{i=4}^{n}a_{i-1}e_i+
\sum\limits_{i=4}^{n}\beta_{i-1}e_i+a_2\gamma
e_n+\sum\limits_{i=3}^{n-2}a_i\sum\limits_{j=i+2}^{n}
\beta_{j-i+1}e_j,$

$[y_n,y_n]=b_1^{2}e_3+b_1\sum\limits_{i=4}^{n}b_{i-1}e_i+
b_1\sum\limits_{i=4}^{n}\beta_{i-1}e_i+\gamma
e_n+\sum\limits_{i=3}^{n-2}b_i\sum\limits_{j=i+2}^{n}
\beta_{j-i+1}e_j.$

If $b_1\neq0$ then $p=1$ and the gradation of the algebra is
natural, i.e we obtain the algebra $NGF_2$. Further we may suppose
that $b_1=0.$

Thus, using lemma \ref{lem1} we have

$[y_n,y_1]=a_2\gamma e_n+\sum\limits_{i=4}^{n}b_{i-1}e_i+
a_2\sum\limits_{i=3}^{n-2}b_i\sum\limits_{j=i+2}^{n}
\beta_{j-i+1}e_j=a_2\gamma e_n+b_3e_4+\\
\sum\limits_{j=5}^{n}(b_{j-1}+
a_2\sum\limits_{i=3}^{j-2}b_i\beta_{j-i+1})e_j,$

$[y_1,y_n]=a_2\gamma e_n+\sum\limits_{i=4}^{n}\beta_{i-1}e_i+
\sum\limits_{i=3}^{n-2}a_i\sum\limits_{j=i+2}^{n}
\beta_{j-i+1}e_j=a_2\gamma e_n+\beta_3e_4+\\
\sum\limits_{j=5}^{n}(\beta_{j-1}+
\sum\limits_{i=3}^{j-2}a_i\beta_{j-i+1})e_j,$

$[y_n,y_n]=\gamma e_n+\sum\limits_{i=3}^{n-2}b_i
\sum\limits_{j=i+2}^{n}\beta_{j-i+1}e_j= \gamma e_n+
\sum\limits_{j=5}^{n}(\sum\limits_{i=3}^{j-2}b_i\beta_{j-i+1})e_j.$

\textbf{Case 1.} There exists some $i$ $(3\leq i \leq n-2)$ such
that $\beta_i\neq0.$ Let $\beta_k$ be the first non zero coefficient
from the set $\{\beta_3, \beta_4, \dots, \beta_{n-2}\}.$

Then

$[y_n,y_1]=a_2\gamma e_n+\sum\limits_{j=4}^{k+1}b_{j-1}e_j+
\sum\limits_{j=k+2}^{n}(b_{j-1}+a_2\sum\limits_{i=3}^{j+1-k}
b_i\beta_{j-i+1})e_j,$

$[y_1,y_n]=a_2\gamma e_n+\beta_{k}e_{k+1}+\sum\limits_{j=k+2}^{n}
(\beta_{j-1}+\sum\limits_{i=3}^{j+1-k}a_i\beta_{j-i+1})e_{j},$

$[y_n,y_n]=\gamma e_n+
\sum\limits_{j=k+2}^{n}(\sum\limits_{i=3}^{j+1-k}b_i\beta_{j-i+1})e_j.$

Since $<y_p, y_n>=V_p,$ from the condition $[V_1,V_p]\subseteq
V_{p+1}=<y_{p+1}>$ it follows that

$[y_1,y_n]=\delta y_{p+1}=\delta(e_{p+2}+(*)e_{p+3}+\dots+(*)e_n).$

On the other hand

$[y_1,y_n]=a_2\gamma e_n+\beta_{k}e_{k+1}+ \sum\limits_{j=k+2}^{n}
(\beta_{j-1}+\sum\limits_{i=3}^{j+1-k}a_i\beta_{j-i+1})e_{j}.$

Comparing the coefficients we obtain $\delta=\beta_k$ and $p=k-1.$

Thus we have $V_i=<y_i> \ {for} \ 1 \leq i \leq n-1, \ i\neq k-1,$
$V_{k-1}=<y_{k-1}, y_n>.$

Therefore

$[y_n,y_1]=a_2\gamma e_n+\sum\limits_{j=4}^{k+1}b_{j-1}e_j+
\sum\limits_{j=k+2}^{n}(b_{j-1}+a_2\sum\limits_{i=3}^{j+1-k}
b_i\beta_{j-i+1})e_j,$

but on the other hand from $[y_n,y_1]\in V_{k},$ it follows
 $$[y_n,y_1]=\lambda y_k=\lambda(e_{k+1}+(*)e_{k+2}+\dots+e_n).$$
Therefore $b_i=0, \ 3 \leq i \leq k-1$ and $\lambda=b_k,$ i.e.
$[y_n,y_1]=b_ky_{k}.$

From the equality
$$[y_1,[y_1,y_n]]=[[y_1,y_1],y_n]-[[y_1,y_n],y_1]$$
and the fact that $[y_1,y_n]\in \mathfrak{R}(L)$ we have
$[[y_1,y_1],y_n]=[[y_1,y_n],y_1].$ Therefore
$[y_2,y_n]=[[y_1,y_n],y_1]=\beta_ky_{k+1}.$

From the equalities
$$[y_i,[y_1,y_n]]=[[y_i,y_1],y_n]-[[y_i,y_n],y_1]$$
by induction we obtain that $[y_i,y_n]=\beta_ky_{k+i-1}, \ 2\leq i
\leq n-k.$

Therefore we have the following multiplications in the algebra $L.$

$$\left\{\begin{array}{ll}
[y_i,y_1]=y_{i+1}, &  \ 1\leq i \leq n-2, \\[1mm]
[y_n,y_1]=b_ky_{k}, & \\[1mm]
[y_i,y_n]=\beta_ky_{k+i-1},& \ 1 \leq i \leq n-k,\\[1mm]
[y_n,y_n]=\delta y_{2k-2}.&
\end{array} \right.$$

\textbf{Case 1.1.} Let $2k-2 < n-1.$ Then
$[y_n,y_n]=b_k\beta_{k}y_{2k-2.}$ Taking a change of the basic
element as follows $y_n^{'}=\frac{1}{\beta_k}(y_n-b_ky_{k-1}),$ we
obtain $b_k=0, \ \beta_k=1.$ Thus we obtain the algebra $M_1(k), \ \
3 \leq k \leq n-2, \ k\neq\frac{n+1}{2}.$

\textbf{Case 1.2.} Let $2k-2=n-1.$ Then
$[y_n,y_n]=(b_k\beta_{k}+\gamma)y_{n-1}$ and changing the basic
element $y_n$ by $y_n^{'}=\frac{1}{\beta_k}(y_n-b_ky_{k-1}),$ we
obtain $b_k=0, \ \beta_k=1,$ i.e. $[y_n,y_n]=\gamma y_{n-1}.$ If
$\gamma=0,$ we obtain the algebra $M_1(\frac{n+1}{2}).$ If $\gamma
\neq 0,$ then we get the algebra $M_2.$

\textbf{Case 2.} Let $\beta_i=0$ for any $i$ ($3\leq i\leq n-2.$)
Then

$[y_n, y_1]=a_2\gamma e_n + \sum\limits_{i=4}^{n}b_{i-1}e_i,$

$[y_1, y_n]=(a_2\gamma + \beta_{n-1})e_n,$

$[y_n, y_n]=\gamma e_n.$

\textbf{Case 2.1.} There exists $i$ ($3\leq i \leq n-2$) such that
$b_i\neq 0.$ Let $b_k$ be the first non zero element from the set
$\{b_3, b_4, \dots, b_{n-2}\}.$ Then by the same arguments as in the
case 1 we obtain an algebra with the following table of
multiplication:

$$\left\{\begin{array}{ll}
[y_i, y_1]=y_{i+1}, &  \ 1 \leq i \leq n-2, \\[1mm]
[y_n,y_1]=b_ky_{k}, & \\[1mm]
[y_i, y_n]=0, & \ 1 \leq i \leq n-1,\\[1mm]
[y_n, y_n]=\gamma y_{n-1}. &
\end{array} \right.$$

It is not difficult to see that taking $y_n^{'}=y_n-b_ky_{k-1}$ we
obtain $b_k=0.$ In the case $\gamma=0$ we have the algebra $NGF_2.$
If $\gamma\neq0,$ then by scaling the basis we obtain $\gamma=1,$
i.e. the algebra $M_3.$

\textbf{Case 2.2.}  Let $b_i=0$ for any $i$ ($3\leq i \leq n-2.$)
Then we have

$[y_n, y_1]=(a_2\gamma + b_{n-1})e_n,$

$[y_1, y_n]=(a_2\gamma+\beta_{n-1})e_n,$

$[y_n, y_n]=\gamma e_n.$

Suppose that $\gamma=0,$ then in the case
($b_{n-1},\beta_{n-1})=(0,0)$ we have the algebra $NGF_2$. If
($b_{n-1},\beta_{n-1})\neq(0,0)$ then by a change of basis (similar
to the case 2.1.) we obtain either the algebra $NGF_2$ or
$M_1(n-1).$

Let $\gamma\neq 0,$ then $n$ is odd and $p=\frac{n-1}{2}.$ If
$a_2\gamma+\beta_{n-1}\neq 0$ or $a_2\gamma+\b_{n-1}\neq 0$ then
from the conditions $[y_1,y_n], \ [y_n,y_n]\in V_{n-1}$ we have
$p=1, n=3,$ i.e. a three-dimensional naturally graded algebra. So,
we have to consider the case
 $a_2\gamma+\beta_{n-1}=a_2\gamma+\b_{n-1}=0.$ Taking the following change of
 basis:
  $y_n^{'}=\frac{y_1}{\sqrt{\gamma}}$ we may suppose that $\gamma=1,$
i.e. we obtain the algebra $M_3.$ \hfill{\fbox{}}

From \cite{CGO} we know that the family $F_2$ of theorem \ref{th2}
contains no algebra which admits a gradation of length $n.$
Therefore in the proposition \ref{prop10} we obtain the
classification of filiform algebras of length $(n-1)$ from the
family $F_2.$

\begin{prop} \label{prop3}

Let $L$ be any filiform Leibniz algebra from the family $F_3$ in
 theorem \ref{th2}. If $L$ admits a gradation with length
$(n-1),$ then it is isomorphic either to the algebra
$$NGF_3=\left\{\begin{array}{lll} [y_i,y_1]=-[y_1,y_i]=y_{i+1}, &
2\leq i \leq {n-1}\\[1mm]
[y_i,y_{n+1-i}]=-[y_{n+1-i},y_i]=\alpha (-1)^{i+1}y_n & 2\leq i\leq
n-1.
\end{array} \right.$$ or to the algebra
$$M_4=\left\{\begin{array}{ll}
[y_1,y_1]=y_{n}, & \\[1mm]
[y_i,y_1]=y_{i+1}, &  \  2\leq i \leq n-1, \\[1mm]
[y_1,y_i]=-y_{i+1}, & \  2\leq i \leq n-1.
\end{array} \right.$$
where all omitted products are equal to zero, $\alpha\in\{0,1\}$ for
even $n$ and $\alpha=0$ for odd $n,$ and $\{y_1, y_2, \ \dots,
y_n\}$ is a basis of the corresponding algebra.
\end{prop} \noindent {\bf Proof.} As in the above propositions we choose
generators from the homogeneous basis

$x_s=e_1+a_2e_2+ \dots + a_ne_n$ and $x_t=b_1e_1+e_2+b_3e_3+ \dots
+b_ne_n.$

Consider the products
$$[x_s,x_s]=(\theta_1+a_2\theta_2+a_2^{2}\theta_3)e_n,$$
$$[x_t,x_t]=(b_1^{2}\theta_1+b_1\theta_2+\theta_3)e_n,$$
$$[\dots[x_t,\underbrace{x_s], \dots,
x_s]}\limits_{i-times}=(1-a_2b_1)e_{i+2}+(*)e_{i+3}+\dots+(*)e_n, \
1\leq i \leq n-2.$$

Since $1-a_2b_1\neq 0,$ we can choose
$$y_1=x_s, \ y_2=x_t, \ y_i=[\dots[x_t,\underbrace{x_s], \dots,
x_s]}\limits_{(i-2)-times}, \ 3\leq i \leq n.$$ Therefore we have

$<y_1>\subseteq V_{k_s}, \ <y_2>\subseteq V_{k_t}, \ <y_i>\subseteq
V_{k_t+(i-2)k_s} \ 3 \leq i \leq n.$

From the connectedness of the gradation we have that $k_s=\pm 1.$
Moreover, as in the previous cases we can assume that $k_s=1.$ So,
we have $[y_i, y_1]=y_{i+1}, \ 2 \leq i \leq n-1.$

Consider the products

$[y_3, y_2]=b_1(1-a_2b_1)e_4+(*)e_5+\dots +(*)e_n.$

If $b_1\neq0,$ then $[y_3,y_2]=b_1y_4$ and therefore $k_t=1$, i.e.
the gradation coincides with the natural gradation, and we obtain
the algebra $NGF_3.$ Hence we may assume that $b_1=0.$

 Since $l(\oplus L)=n-1,$ there exists $p$ ($1\leq p\leq n-2$) such that
 $V_1=V_{k_{t}+p},$ therefore $k_t=1-p.$

Thus $$L=V_{1-p}\oplus V_{2-p}\oplus \dots \oplus V_{n-1-p},$$ where
$V_{1-p}=<y_2>, \ V_{2-p}=<y_3>, \dots, \ V_0=<y_{p+1}>, \ V_1=<y_1,
y_{p+2}>, \ V_2=<y_{p+3}>, \dots, \ V_{n-1-p}=<y_n>.$

Therefore the table of multiplication in the algebra $L$ is as
follows:
$$\left\{\begin{array}{ll}
[y_1,y_1]=\alpha y_{p+3}, & \\[1mm]
[y_i,y_1]=y_{i+1}, &  \  2\leq i \leq n-1, \\[1mm]
[y_1,y_i]=-y_{i+1}, & \  2\leq i \leq n-1, \\[1mm]
[y_i,y_j]=0 & \ 2\leq i \leq n, \ 2\leq j \leq n.
\end{array} \right.$$

In the case where $\alpha=0$ without loss of generality we may
 suppose that $p=0$ and therefore we have the natural
 gradation and thus we obtain the algebra $NGF_3.$

If $\alpha\neq0$ then $p=n-3$ and taking the change of basis
$$y_1^{'}=y_1, \ y_i^{'}=\alpha y_i, \ 2\leq i \leq n$$ we can
assume that $\alpha=1,$ i.e. we obtain the algebra
$M_4.$\hfill{\fbox{}}

From \cite{GJR1} it follows that the algebra $NGF_3$ only for
$\alpha=0$ admits a gradation with the length $n,$ i.e. this has the
maximal length. From \cite{CGO} we have that $M_4$ does not admit
any gradation with the length $n.$ Therefore the algebras $NGF_3$
for $\alpha=1$ and $M_4$ are algebras of length $(n-1).$

Summarizing the propositions \ref{prop1} - \ref{prop3} we obtain the
following result.

\begin{thm} \label{th3} Any $n$-dimensional complex filiform
Leibniz algebra of length $(n-1)$ is isomorphic to one of the
following pairwise non isomorphic algebras:
$$NGF_2, \ NGF_3(\alpha=1) , \ M_1(k), M_2, \ M_3, \ M_4.$$
\end{thm}
\noindent {\bf Proof.} The algebra $M_2$ is not isomorphic to the
algebra $M_1(k)$ for any $k$ ($3\leq k \leq n-2$), because
$dim\mathfrak{L}(M_1(k))=2$ but $dim\mathfrak{L}(M_2)=1.$ The fact
that the algebras $M_1(k), \ M_2, \ M_3, \ M_4$ are pairwise non
isomorphic follows from the criteria of isomorphism of two filiform
algebras of the class $F_2$ (\cite{GO}, theorem 4.4). Since the
families of algebras $F_1, \ F_2, \ F_3$ are disjoint, $NGF_2$ is a
split non Lie Leibniz algebra, and $NGF_3(\alpha=1)$ is a Lie
algebra, the proof is complete. \hfill{\fbox{}}

\section{Derivations of $n$-dimensional filiform Leibniz
\\ algebras of the length $(n-1)$.}

In this section we apply the above classification of $n$-dimensional
filiform Leibniz algebras of length $(n-1)$ to the description of
their derivations.

Let $L$ be a Z-graded Leibniz algebra, i.e.
$L=\bigoplus\limits_{i\in Z}V_i.$ This gradation induces a gradation
of the algebra $Der(L)=\bigoplus\limits_{i\in Z}W_i$ in the
following way:

$$W_i=\lbrace d\in Der(L) \mid d(x)\in W_{i+j} \ {for \  any} \ x\in V_j\rbrace.$$

If the gradation of $L$ is finite, then the gradation of $Der(L)$ is
also finite. In particular, if for an $n-$dimensional filiform
Leibniz algebra $L$ we have $L=V_{k_1}\oplus V_{k_1+1}\dots\oplus
V_{k_1+n-2}$ then it is easy to see that $Der(L)=W_{2-n}\oplus
W_{3-n}\dots\oplus W_{n-3}\oplus W_{n-2}.$

Since derivations of the algebras $NGF_1-NGF_3$ have already been
described in
 \cite{O}, \cite{GK} we only have to consider the remaining algebras from  theorem \ref{th3}.

First, consider the family of algebras

$$M_1(k)=\left \{\begin {array}{ll}
[y_i,y_1]=y_{i+1}, &  \  1\leq i \leq n-2 \\[1mm]
[y_i,y_n]=y_{k+i-1}, & \  1\leq i \leq n-k, \  3\leq k \leq n-2,
\end{array} \right.$$

\begin{prop} \label{prop5}
The linear maps $h_i, \ d_j\ (1\leq i \leq 2, \ 0\leq j \leq n-2)$
on $M_1(K)$ defined as follows:

$ d_0(y_i)=iy_i, \ 1 \leq i \leq n-1, \ d_0(y_n)=(k-1)y_n,$

$d_j(y_i)=y_{j+i}, \ 1 \leq j \leq n-2, \ 1 \leq i \leq n-j-1,$

$h_1(y_n)=y_{n-1},$

in the case $2k-2\leq n-1$ \\ and as

$ d_0(y_i)=iy_i, \ 1 \leq i \leq n-1, \ d_0(y_n)=(k-1)y_n,$

$d_j(y_i)=y_{i+j}, \ 1 \leq j \leq n-2, \ 1 \leq i \leq n-j-1,$

$h_1(y_n)=y_{n-1},$

$h_2(y_1)=y_{n}, \ h_2(y_i)=(i-1)y_{k+i-2}, \ 2 \leq i \leq n-k+1,$

in the case $2k-2 \geq n$ \\
form a basis of the space $Der(M_1(k)).$

\end{prop}

\noindent {\bf Proof.} It is easy to see that the given maps are
derivations and they are linearly independent.

From proposition \ref{prop10} we have $M_1(K)=V_1\oplus
V_2\oplus\dots \oplus V_{n-1},$ where $V_i=<y_i> \ 1 \leq i \leq
n-1, i\neq k-1$ and $V_{k-1}=<y_{k-1}, y_n>, \ 3 \leq k \leq n-1.$

Let $d\in Der(M_1(k)),$ then we have the decomposition
$d=\sum\limits_{i=2-n}^{n-2}d_i.$

It is clear that $d_{2-n}(y_i)=0$ for $1\leq i \leq n-2,$
$d_{2-n}(y_n)=0$ and $d_{2-n}(y_{n-1})=\alpha y_1, \alpha\in C.$

The property of derivations implies

$d_{2-n}(y_{n-1})=d_{2-n}([y_{n-2},y_1])=
[d_{2-n}(y_{n-2}),y_1]+[y_{n-2},d_{2-n}(y_1)]=0,$ \\
i.e. $d_{2-n}=0.$

Similarly we obtain that $d_j=0, \ j\leq 1-k.$

Consider $d_{2-k}\in Der(M_1(k)).$ It is easy to see that
 $d_{2-k}(y_i)=0$ for $1\leq i \leq k-2,$ and $d_{2-k}(y_n)=\alpha
y_1$ for some $\alpha\in C.$

From

$d_{2-k}(y_i)=d_{2-k}([y_{i-1},y_1])=[d_{2-k}(y_{i-1}),y_1]+
[y_{i-1},d_{2-k}(y_{1})]=0,  \ 2 \leq i \leq n-1,$

we have $d_{2-k}(y_i)=0, \ 2 \leq i \leq n-1.$

By the property of derivations one has

$0=d_{2-k}([y_n,y_1])=[d_{2-k}(y_n),y_1]+[y_n,d_{2-k}(y_1)]= [\alpha
y_1,y_1]=\alpha y_2 \Rightarrow \alpha=0,$ i.e. $d_{2-k}=0.$

Analogously we obtain $d_{t-k}=0, \ 2 \leq t \leq k-1.$ Thus we have
$d_i=0, \ 2-n \leq i \leq -1.$

Consider $d_0\in Der(M_1(k)).$ Then
$$d_0(y_i)=\left\{\begin{array}{ll}
\beta_{0, i}y_i, &  \  1\leq i \leq n-1, \  i\neq k-1, \\[1mm]
\gamma_{0, 1}y_n+\gamma_{0, 2}y_{k-1}, & \  i=n,\\[1mm]
\gamma_{0, 3}y_n+\gamma_{0, 4}y_{k-1}, & \  i=k-1,
\end{array} \right.$$
where $\beta_{0, i},\gamma_{0, j}$ are complex numbers.

Further we have

$d_0(y_2)=d_0([y_1, y_1])=[d_0(y_1),y_1]+[y_1,d_0(y_1)]=[\beta_{0, 1}y_1, y_1]+[y_1,\beta_{0,1}y_1]=2\beta_{0, 1}y_2 \Rightarrow \beta_{0,2}=2\beta_{0,1}.$

From the chain of equalities
$$d_0(y_i)=d_0([y_{i-1},y_1])=[d_0(y_{i-1}),y_1]+[y_{i-1},d_0(y_1)]=$$
$$[\beta_{o,
i-1}y_{i-1},y_1]+[y_{i-1},\beta_{0,1}y_1]=
(\beta_{0,i-1}+\beta_{0,1})=i\beta_{0,1}y_i,$$ by induction we
obtain $\beta_{0,i}=i\beta_{0,1}, \ 2 \leq i \leq n-1, \ i\neq k-1$
and $\gamma_{0,3}=0, \ \gamma_{0,4}=(k-1)\beta_{0,1}.$

The property of derivations implies

$0=d_0([y_n,y_1])=[d_0(y_n),y_1]+[y_n,d_0(y_1)]=[\gamma_{0,1}y_n+\gamma_{0,2}y_{k-1},y_1]+[y_n,\beta_{0,1}y_1]=
\gamma_{0,2}y_k \Rightarrow \gamma_{0,2}=0.$

Further we have
$$d_0([y_1,y_n])=[d_0(y_1),y_n]+[y_1,d_0(y_n)]=[\beta_{0,1}y_1,y_n]
+[y_1,\gamma_{0,1}y_n]=(\beta_{0,1}+\gamma_{0,1})y_k.$$ On the other
hand
$$d_0([y_1,y_n])=d_0(y_k)=k\beta_{0,1}y_k.$$
Therefore $\gamma_{0,1}=(k-1)\beta_{0,1}.$

Thus we have
$$d_0(y_i)=\left\{\begin{array}{ll}
iy_i, &  \  1\leq i \leq n-1, \\[1mm]
(k-1)y_n, & \  i=n.
\end{array} \right.$$

Consider $d_j\in Der(M_1(k)), \ 1 \leq j \leq n-2.$ It is clear that

$$d_j(y_i)=\left\{\begin{array}{ll}
\beta_{j, i}y_{i+j}, &  \  1\leq i \leq n-1, \  i+j\neq k-1, \\[1mm]
\gamma_{j, 1}y_n+\gamma_{j, 2}y_{k-1}, & \  i+j=k-1, \ \ \quad \quad \quad \quad \quad \quad \quad \quad \quad \quad \quad \quad (1) \\[1mm]
\gamma_{j, 3}y_{k-1+j} & \  i=n,
\end{array} \right.$$

We have

$0=d_j([y_n, y_1])=[d_j(y_n), y_1]+[y_n, d_j(y_1)]=[\gamma_{j,3}y_{k-1+j},y_1]+[y_n, d_j(y_1)]=\gamma_{j,3}y_{k+j}.$

If $j\leq n-k-1,$ then $\gamma_{j,3}=0$ and if $j=n-k,$ then
$\gamma_{j,3}$ is arbitrary. In the case where $j\geq n-k+1$ one has
$d_j(y_n)=0.$

Consider the derivation $d_{k-2}.$ From (1) we have

$d_{k-2}(y_1)=\gamma_{k-2, 1}y_{n}+\gamma_{k-2, 2}y_{k-1},$

$d_{k-2}(y_2)=d_{k-2}([y_1, y_1])=
[d_{k-2}(y_1),y_1]+[y_1,d_{k-2}(y_1)]=
[\gamma_{k-2,1}y_{n}+\gamma_{k-2,2}y_{k-1},y_1]+
[y_1,\gamma_{k-2,1}y_{n}+\gamma_{k-2,2}y_{k-1}]=(\gamma_{k-2,1}+\gamma_{k-2,2})y_k,$

$d_{k-2}(y_3)=d_{k-2}([y_2, y_1])=
[d_{k-2}(y_2),y_1]+[y_2,d_{k-2}(y_1)]=
[(\gamma_{k-2,1}+\gamma_{k-2,2})y_k,y_1]+ [y_2,
\gamma_{k-2,1}y_{n}+\gamma_{k-2,2}y_{k-1}]=
(2\gamma_{k-2,1}+\gamma_{k-2,2})y_{k+1}.$

By induction we obtain
$$d_{k-2}(y_i)=d_{k-2}([y_{i-1},y_1])=
[d_{k-2}(y_{i-1}),y_1]+[y_{i-1}, d_{k-2}(y_1)]=$$
$$[((i-2)\gamma_{k-2,1}+\gamma_{k-2,2})y_{k+i-3},y_1]+
[y_{i-1},\gamma_{k-2,1}y_{n}+\gamma_{k-2,2}y_{k-1}]=$$
$$((i-1)\gamma_{k-2,1}+\gamma_{k-2,2})y_{k+i-2}, \ 2 \leq i \leq
n-k+1. \quad  \quad  \quad \quad \quad (2)$$

On the other hand
$$d_{k-2}(y_k)=d_{k-2}([y_1,y_n])=
[d_{k-2}(y_1),y_n]+[y_1,d_{k-2}(y_n)]=$$
$$[\gamma_{k-2,1}y_{n}+\gamma_{k-2,2}y_{k-1},y_1]+[y_1,\gamma_{k-2,3}y_{2k-3}]=
\gamma_{k-2,2}y_{2k-2}. \quad \quad \quad (3)
$$
Comparing the coefficients in (2) for $i=k$ with the coefficients in
(3) we obtain $(k-1)\gamma_{k-2,1}+\gamma_{k-2,2}=\gamma_{k-2,1}$
for $1 \leq 2k-2\leq n-1,$ therefore $\gamma_{k-2,1}=0.$ For $2k-2 >
n-1,$ we have
$d_{k-2}(y_i)=((i-1)\gamma_{k-2,1}+\gamma_{k-2,2})y_{k+i-2}, \ 2\leq
i \leq n-k+1.$

\textbf{Case 1.} $2k-2\leq n-1.$ In this case
 $d_{k-2}(y_i)=\gamma_{k-2,2}y_{k+i-2}, \ 1 \leq i \leq n-k+1.$

Consider $d_j(y_1)=\beta_{j,1}y_{j+1}, \ 2 \leq j \leq n-k+1, \
j\neq k-2.$

Then

$d_{j}(y_2)=d_j([y_1,y_1])=[d_j(y_1),y_1]+[y_1,d_j(y_1)]=
[\beta_{j,1}y_{j+1},y_1]+[y_1,\beta_{j,1}y_{j+1}]=\beta_{j,1}y_{j+2}.$

By induction we obtain

$d_j(y_i)=d_j([y_{i-1},y_1])=[d_j(y_{i-1}),y_1)]+[y_{i-1},d_j(y_1)]=
[\beta_{j,1}y_{j+i-1},y_1]+[y_{i-1},\beta_{j,1}y_{j+1}]=\beta_{j,1}y_{j+i},$
i.e. $d_j(y_i)=\beta_{j,1}y_{j+i}, \ 1 \leq j \leq n-2, \ 1 \leq i \leq
n-1-j, j \neq n-k.$

By the above consideration we can write $d_{n-k}$ in the following
form $d_{n-k}=\beta_{n-k,1}d_{n-k}^{'}+\gamma_{n-k,3}h_1,$ where
$d_{n-k}^{'}(y_i)=y_{n-k+i}, \ h_1(y_n)=y_{n-1}.$

Denoting $d_{n-k}^{'}$ by $d_{n-k}$ we obtain $d=h_1+d_0+d_1+\dots +
d_{n-2}.$

\textbf{Case 2.} Let $2k-2 \geq n,$ then

$d_{k-2}(y_1)=\gamma_{k-2,1}y_{n}+\gamma_{k-2,2}y_{k-1},$

$d_{k-2}(y_i)=((i-1)\gamma_{k-2,1}+\gamma_{k-2,2})y_{k+i-2}, \ 2
\leq i \leq n-k+1,$ i.e.

$d_{j}(y_i)=\beta_{j,i}y_{j+i}, \ 1\leq j \leq
n-2, \ 1 \leq i \leq n-1-j, \ j \notin\{n-k, k-2\}.$

As above one has
 $d_{k-2}=\gamma_{k-2,2}d_{k-2}^{'}+\gamma_{k-2,1}h_2,$ where
$d_{k-2}^{'}(y_i)=y_{k-2+i}, \ h_2(y_1)=y_n$ and
$h_2(y_i)=(i-1)y_{k+i-2}.$

Similarly $d_{n-k}=\beta_{n-k,1}d_{n-k}^{'}+\gamma_{k-2,3}h_1,$
where $d_{n-k}^{'}(y_i)=y_{n-k+i}, \ h_1(y_n)=y_{n-1}.$

Denoting $d_{k-2}^{'}$ by $d_{k-2}$ and $d_{n-k}^{'}$ by $d_{n-k}$
we obtain that $d=h_1+h_2+d_0+d_1+\dots+d_{n-2}.$ \hfill{\fbox{}}

Similarly one can prove the following propositions.

\begin{prop}

The linear maps $h_1, \ d_j, \ (0 \leq j \leq n-2)$ on $M_2$
(respectively on $M_3$) defined as

$d_0(y_n)=(k-1)y_n, \ d_0(y_i)=iy_i, \ 1 \leq i \leq n-1,$

$d_j(y_i)=y_{i+j}, \ 1 \leq j \leq n-2, \ 1 \leq i \leq n-j-1,$

$h_1(y_n)=y_{n-1},$

form a basis of the space $Der(M_2)$ (respectively $Der(M_3)$).

\end{prop}

\begin{prop}

The linear maps $h_0, \ h_1, \ d_j, \ (3-n \leq j \leq n-2)$ on
 $M_4$ defined as

$d_{j-n}(y_1)=y_{j-1}, \ 3 \leq j \leq n,$

$d_j(y_i)=y_{i+j}, \ 1 \leq j \leq n-2, \ 2 \leq i \leq n-j,$

$h_0(y_1)=y_{n}, \ h_0(y_i)=(2-n+i)y_i, \ 2 \leq i \leq n,$

$h_1(y_1)=y_{n}.$

form a basis of the space $Der(M_4).$
\end{prop}

Recall \cite{L} that the 1-cocycle space $Z^1(L,L)$ of the algebra
$L$ with the values in $L$ is the space $Der(L)$ of derivations of
$L$, while the 1-coboundary space $B^1(L,L)$ is the space $Inn(L)$
of inner derivations. The first cohomology group $H^1(L,L)$ is the
quotient space $Z^1(L,L)/ B^1(L,L).$ Thus the above results imply

\begin{cor} The first cohomology groups of Leibniz algebras of
length $(n-1)$ have the following dimensions:

$DimH^1(M_1(k), M_1(k))=n-2,\ \mbox{for} \ 2k-2\leq n-1,$

$DimH^1(M_1(k), M_1(k))=n-1,\ \mbox{for} \ 2k-2\geq n,$

$DimH^1(M_i, M_i)=n-2,\ \mbox{for} \  2\leq i \leq 3,$

$DimH^1(M_4, M_4)=n-3.$

\end{cor}

\begin{rem} The cases of the algebras $NGF_2$ and $NGF_3\
(\alpha=1)$ have been already considered in \cite{O} and \cite{GK},
respectively.
\end{rem}

Further recall \cite{L} that the 2-coboundary space for the algebra
$L$ is $B^2(L,L)= \{f:L\otimes L\rightarrow L  \mid
f(x,y)=[d(x),y]+[x,d(y)]-d([x,y])\  \mbox{for \ some \ linear \
transformation  \ d  \ of}\ L\}.$  For the algebras $NGF_2$ and
$NGF_3\ (\alpha=1)$ 2-coboundary spaces were considered in \cite{O}
and \cite{GK}. For the rest of Leibniz algebras of length $(n-1)$ we
have

\begin{cor}
$DimB^2(M_1(k), M_1(k))=n^2-n+2,\ \mbox{for} \ 2k-2\leq n-1,$

$DimB^2(M_1(k), M_1(k))=n^2-n+1,\ \mbox{for} \ 2k-2\geq n,$

$DimB^2(M_i, M_i)=n^2-n+2, \ \mbox{for} \  2\leq i \leq 3,$

$DimB^2(M_4, M_4)=n^2-n+3.$

\end{cor}
\textbf{Acknowledgments.} \emph{The second and third named authors
would like to acknowledge the hospitality of the $\,$ "Institut
f\"{u}r Angewandte Mathematik",$\,$ Universit\"{a}t Bonn (Germany).
This work is supported in part by the DFG 436 USB 113/10/0-1 project
(Germany) and the Fundamental Research Foundation of the Uzbekistan
Academy of Sciences.}

\end{document}